\documentclass{amsart}
\usepackage{amscd,amssymb}
\usepackage[mathcal]{euscript}
\newtheorem{theorem}{Theorem}[section]
\newtheorem{lemma}[theorem]{Lemma}
\newtheorem{proposition}[theorem]{Proposition}
\newtheorem{axiom}[theorem]{Axiom}
\newtheorem{corollary}[theorem]{Corollary}
\theoremstyle{definition}
\newtheorem{definition}[theorem]{Definition}
\newtheorem{remark}[theorem]{Remark}

\numberwithin{equation}{subsection}

\DeclareMathOperator{\HR}{\Ho(Reedy)}
\DeclareMathOperator{\HC}{Ho(can.)}
\DeclareMathOperator{\Ho}{Ho}

\DeclareMathOperator{\map}{map}
\DeclareMathOperator{\Sing}{Sing}
\DeclareMathOperator{\hocolim}{hocolim}

\DeclareMathOperator{\colim}{colim}
\DeclareMathOperator{\oplim}{lim}
\newcommand{\parallelarrows}[1]{\begin{array}{c} {\hbox to
#1{\rightarrowfill}}  \vspace{-0.35cm} \\ {\hbox to
#1{\rightarrowfill}} \end{array}}

\newcommand{\del}{\partial}
\newcommand{\boxprod}{\mathbin\square}
\newcommand{\simp}{realization\ }
\newcommand{\real}{Realization\ }

\newcommand{\iso}{\cong}

\newcommand{\Ss}{\mathcal{S}}

\newcommand{\M}{\mathcal{M}}
\newcommand{\cL}{\mathcal{L}}

\newcommand{\mc}{\colon\,}
\newcommand{\mathcolon}{\colon\,}
\newcommand{\s}{\mathcolon}

\newcommand{\C}{\mathcal{C}}
\newcommand{\D}{\mathcal{D}}
\newcommand{\Z}{\mathbb{Z}}
\newcommand{\sD}{s\mathcal{D}}

\newcommand{\id}{{\mbox {id}}}

\newcommand{\ev}{{\mbox {Ev}}}
\newcommand{\osum}{\oplus}
\newcommand{\bU}{\bar{U}}
\newcommand{\mbD}{{\bar \Delta}}

\title[]{Simplicial structures \\on model categories and functors} 

\date{\today}
\author{Charles Rezk}
\thanks{The first author was partially supported by an AMS Centennial 
Fellowship}
\address{Institute for Advanced Study \\ School of Mathematics\\
Olden Lane\\ Princeton, NJ 08540\\ USA}
\email{rezk@ias.edu}
\author{Stefan Schwede}
\address{Fakult\"at f\"ur Mathematik \\ Universit\"at Bielefeld\\
33615 Bielefeld, Germany}
\email{schwede@mathematik.uni-bielefeld.de}
\author{Brooke Shipley}
\thanks{The third author was partially supported by NSF grants} 
\address{Department of Mathematics \\ Purdue University \\ West Lafayette, IN 
47907
\\ USA} 
\email{bshipley@math.purdue.edu}

\begin{document}

\begin{abstract}
We produce a highly structured way of associating a simplicial category to a
model category which improves on work of Dwyer and Kan and answers 
a question of Hovey.  
We show that model categories satisfying a certain axiom are Quillen
equivalent to simplicial model categories.  A simplicial model category
provides higher order structure such as composable mapping spaces 
and homotopy colimits.  We also show that certain
homotopy invariant functors can be replaced by weakly equivalent simplicial,
or `continuous', functors.  This is used to show that if a
simplicial model category structure exists on a model category then it is 
unique up to simplicial Quillen equivalence. 
\end{abstract}

\maketitle

\section{Introduction}\label{intro}
In~\cite{DK} Dwyer and Kan showed that a simplicial category, called the 
hammock localization, can be associated
to any Quillen model category~\cite{Q}.  This simplicial category
captures higher order information, for example fibration and cofibration
sequences and mapping spaces, see~\cite[I 3]{Q}, which is not captured by the 
ordinary homotopy category.  Hovey carried this further by showing that the 
homotopy category of simplicial sets acts on the homotopy category of any model
category~\cite[5.5.3]{hovey}.  Hovey then wondered if in fact every model 
category is
Quillen equivalent to a simplicial model category~\cite[8.9]{hovey}.  
Quillen equivalence is the appropriate notion of equivalence for
model categories, so
this would be the most highly structured way of associating a simplicial 
category to any model category.  The following existence result is proved  
in Theorem~\ref{main-thm}.
\begin{theorem}\label{thm-1}
If $\C$ is a left proper, cofibrantly generated model category that
satisfies Realization Axiom~\ref{real-axiom}, then $\C$ is Quillen
equivalent to a simplicial model category.
\end{theorem}

Throughout this paper we use a slightly stronger notion of
cofibrantly generated model category than is standard; see 
Definition~\ref{def-cof-gen}.
We also have the following uniqueness result, which is proved as
Corollary~\ref{cor-SMC-Qui}.
Assume that $\C$ and $\D$ are model categories which either satisfy the
hypotheses of Theorem~\ref{thm-1} or satisfy the hypotheses of one of
the general localization machines in~\cite{HH} or~\cite{smith-combinatorial},
see also~\cite{dugger}.
\begin{theorem}\label{thm-2}
Under these hypotheses, if $\C$ and $\D$ are Quillen equivalent simplicial
model categories, then $\C$ and $\D$ are simplicially Quillen equivalent.
\end{theorem}
\noindent By considering the identity functor, this shows
that a simplicial model category structure on a model category is unique up to
simplicial Quillen equivalence, see Corollary~\ref{cor-simp-Qui}.
This strengthens Dwyer and Kan's analogous result on the homotopy categories
in~\cite{DK}.

To prove Theorem~\ref{thm-2}, in Section~\ref{sec-functors} we consider 
replacing functors
between simplicial model categories by simplicial, or `continuous', functors.  
We show that a homotopy invariant functor $F$ can be replaced by a naturally 
weakly equivalent simplicial functor, see Corollary~\ref{all-wk-eq}.
We also show that Quillen adjunctions between simplicial model categories,
the appropriate notion of functors between model categories, can
be replaced by simplicial Quillen adjunctions, see 
Proposition~\ref{prop-Qui-adjoint}.
This answers another part of Hovey's problem,~\cite[8.9]{hovey}.

Another reason to construct replacement simplicial model categories is
to have a simple definition of a homotopy colimit.  The original definition
in~\cite[XII]{BK} generalizes to define a homotopy colimit in any
simplicial model category, see~\cite[20]{HH}.  So the simplicial replacements
considered here provide new
situations where a simple homotopy colimit can be defined.  The Bousfield-Kan
type homotopy colimit on the replacement simplicial model category can be
transported to the original model category via the Quillen equivalence.

Showing that stable model categories have simplicial replacements
was the original motivation for this work, see Section~\ref{examples}. 
\begin{proposition}  
Any proper, cofibrantly generated, stable model 
category is Quillen equivalent to a simplicial model category.  
\end{proposition}

\noindent The category of unbounded differential graded modules over a 
differential graded algebra is one particular
example of a stable model category that was not previously known to have
a Quillen equivalent simplicial replacement.  This example is treated
explicitly in Corollary~\ref{cor-chain complexes} and answers another question 
of Hovey,~\cite[8.9]{hovey}.

For a model category $\C$, our candidate for a Quillen equivalent simplicial
model category is based on the category of simplicial objects in $\C$, $s\C$. 
Reedy~\cite{reedy} establishes the Reedy model category on $s\C$, but
it is neither simplicial nor Quillen equivalent to $\C$, see~\cite[2.6]{e2} or
Corollary~\ref{cor-Reedy-simp}.
So we localize the Reedy model category to create the realization
model category.  Instead of using general machinery to produce the
localization model category, we explicitly define the cofibrations, weak 
equivalences, and fibrations and then check that they form a model category.
This avoids unnecessary hypotheses. In Theorem~\ref{main-thm} we show that
if $\C$ is a left proper, cofibrantly generated model category
that satisfies Realization Axiom~\ref{real-axiom}, then the
realization structure on $s\C$ is a simplicial model category
that is Quillen equivalent to the original model category $\C$.

More generally, we show that there is at most one 
model category on $s\C$ that satisfies certain properties, see 
Theorem~\ref{thm-A}.
When this model category exists on $s\C$  it is Quillen equivalent to
the original model category $\C$, and we refer to it as the
canonical model category structure on $s\C$.  
If $\C$ satisfies the hypotheses of Theorem~\ref{main-thm} as listed
above, then the canonical model category structure on $s\C$ exists
and is simplicial 
since it agrees with the realization model category.  The
applications in Sections~\ref{sec-unique} and~\ref{sec-functors} rely only on 
the existence of the canonical model category on $s\C$ and the fact
that it is simplicial.  

In~\cite{dugger}, Dugger has also developed a way to produce 
replacement simplicial model categories. 
His approach is similar to ours, but he uses the 
two general localization machines that exist for
left proper, cellular model categories, see~\cite{HH} and for
left proper, cofibrantly generated, combinatorial model categories, 
see~\cite{smith-combinatorial}.   
Hence, these hypotheses also ensure the existence of the simplicial, canonical 
model category on $s\C$.  So the applications of Sections~\ref{sec-unique} and 
~\ref{sec-functors} also apply under the conditions investigated 
in~\cite{dugger}.

One drawback with these general machines is that the fibrations 
cannot always be identified in concrete terms.
Our approach here is to explicitly define the fibrations and then
verify the model category axioms.   This approach requires a slightly stronger
notion of ``cofibrantly generated'', see Definition~\ref{def-cof-gen}.
Then for left proper, cofibrantly generated model categories,  
Realization Axiom~\ref{real-axiom} is equivalent
to having the explicit definition of the fibrations, see 
Proposition~\ref{prop-fib-real}.  

{\em Organization:}  
In Section~\ref{begin} we recall the simplicial
structure on $s\C$ and the Reedy model category structure
on $s\C$.  In Section~\ref{main}, we define the canonical model 
category structure on $s\C$, the realization model category structure
on $s\C$, and state the main theorems.  
In Section~\ref{examples} we consider examples including
simplicial model categories, stable model categories, and 
unbounded differential graded modules over a differential graded
algebra.  In Sections~\ref{sec-unique} and~\ref{sec-functors} we consider
the applications mentioned above: the uniqueness of simplicial model
category structures and replacing functors by simplicial functors.
In Section~\ref{sec-Reedy}, we show that the Reedy model
category structure only partially satisfies the compatibility axiom SM7.  
This also gives several statements that are needed in later proofs.  
In Section~\ref{proofs} we verify the main theorem, Theorem~\ref{main-thm},  
which states that the \simp structure on $s\C$ is a simplicial model category
that is Quillen equivalent to the original model category, $\C$.

\section{The Reedy model category for simplicial objects in $\C$}\label{begin}
Here we define the canonical simplicial structure 
on the category of simplicial objects of $\C$, $s\C$.  
This is our candidate category for replacing $\C$
by a simplicial model category.  We also recall the definition of a simplicial 
model category and the Reedy model category structure on $s\C$. 

Let $s\C$ denote the simplicial objects in $\C$, i.e.~the functors 
$\Delta^{op} \to \C$.  Let $\Ss$ denote the category of simplicial sets.
For any category $\C$ with small limits and colimits, $s\C$ is tensored
and cotensored over $\Ss$, compare~\cite[II 1]{Q}.
For a set $S$ and $X \in \C$, let $X \cdot S= \coprod_{s \in S} X.$  
For $X$ in $s\C$ and $K$ in $\Ss$ define $X \otimes K$ in $s\C$ as the
simplicial object with $n$th simplicial degree $(X \otimes K)_n=X_n \cdot K_n$.
For $A$ in $\C$ denote $cA \otimes K$ as $A \otimes K$ in $s\C$ 
where $c \s \C \to s\C$ is the constant object functor.  
Note $cA = A \otimes \Delta[0]$.
The cotensor $X^K$ in $s\C$ is also defined in~\cite[II 1]{Q}.  In this paper 
we mainly use the degree zero part in $\C$ of this cotensor, and denote it 
$X^K$.    
From this simplicial tensor one can define simplicial mapping spaces, 
$\map(X,Y)$ in $\Ss$ for $X, Y \in s\C$ with $n$th simplicial degree 
$\map(X,Y)_n= s\C (X \otimes \Delta[n], Y)$.  So $s\C$ is also enriched
over $\Ss$.

We now recall the definition of a simplicial model category, which
asks that the simplicial structure is compatible with the model category
structure.

\begin{definition} 
A {\em simplicial model category} is a model category $\C$ 
that is enriched, cotensored and tensored over $\Ss$ and
satisfies the following axiom:

\begin{axiom}\cite[II.2 SM7]{Q}\label{SM7}
If $f \mathcolon A \to B$ is a cofibration in $\C$ and $i \mathcolon K \to L$
is a cofibration in $\Ss$ then
\[ q\mathcolon A \otimes L  \coprod_{A\otimes K} B \otimes K  \to B \otimes L \]
\begin{enumerate}
\renewcommand{\labelenumi}{(\arabic{enumi})}
\item is a cofibration;
\item if $f$ is a weak equivalence, then so is $q$;
\item if $i$ is a weak equivalence, then so is $q$.
\end{enumerate}
\end{axiom}
\end{definition}

The first model category we consider on $s\C$ is the Reedy model category 
structure, see~\cite[Theorem A]{reedy} or~\cite[2.4]{e2}.
Before defining the Reedy model category structure we need to define
latching and matching objects.
Let $\cL_n$ be the category with objects the maps $[j] \to [n] \in \Delta^{op}$
with $j < n$ and with morphisms the commuting triangles.  Let $l\mathcolon
\cL_n \to \Delta^{op}$ be the forgetful functor.   Given $X \mathcolon
\Delta^{op} \to  \C$, an object in $s\C$, define $L_n X = \colim_{\cL_n} l^* X$.
$L_n X$ is the $n$th {\em latching object} of $X$.
Similarly, let $\M_n$ be the category with objects the maps $[n] \to [j] \in
\Delta^{op}$
with $j < n$ and with morphisms the commuting triangles.  Let $m\mathcolon
\M_n \to \Delta^{op}$ be the forgetful functor.   Given $X \mathcolon
\Delta^{op} \to  \C$, an object in $s\C$, define $M_n X = \oplim_{\M_n} m^* X$.
$M_n X$ is the $n$th {\em matching object} of $X$.

\begin{definition}\label{def-Reedy}
A map $f \mc X \to Y$ in $s\C$ is a {\em level weak equivalence} 
if $X_n \to Y_n$
is a weak equivalence in $\C$ for each $n$.  It is a {\em Reedy cofibration} if
the induced map $X_n \coprod_{L_n X} L_n Y \to Y_n$ is a cofibration
in $\C$ for each $n$.  Similarly, $f$ is a {\em Reedy fibration} 
if the induced map $X_n \to Y_n \prod_{M_n Y} M_n X$ is a fibration in $\C$.
\end{definition}

Note that a map $X \to Y$ in $s\C$ is a
Reedy trivial cofibration (resp.\ Reedy trivial fibration) if and only if
all the maps  $X_n \coprod_{L_n X} L_n Y \to Y_n$ are acyclic cofibrations
in $\C$ (resp.\ all the maps $X_n \to Y_n \prod_{M_n Y} M_n X$ 
are acyclic fibrations in $\C$).
The following theorem is due to Reedy,~\cite[Theorem A]{reedy}.  
See also~\cite[2.4]{e2} or~\cite[5.2.5]{hovey}.

\begin{theorem}\label{thm-Reedy}
The category $s\C$ equipped with the level weak equivalences, 
Reedy cofibrations, and Reedy fibrations is a model category, 
referred to as the {\em Reedy model category}.
\end{theorem}

This Reedy model category structure on $s\C$ with the canonical simplicial 
structure described above 
satisfies properties $(1)$ and $(2)$ of Axiom~\ref{SM7} (SM7) but does not 
satisfy property $(3)$.  This is stated in Corollary~\ref{cor-Reedy-simp}. 
So this model category is not a simplicial model category, 
but is a stepping stone for
defining the model category structure on $s\C$ that is simplicial.

\section{Statement of results}\label{main}
Here we define the realization model category structure on $s\C$.
This is the model category structure on $s\C$ which is simplicial 
and also Quillen equivalent to the original model category on $\C$, see 
Theorem~\ref{main-thm}.
We first show that there is at most one model category on $s\C$ 
with certain properties, which we call the canonical model category, see 
Theorem~\ref{thm-A}.  We then show that the 
canonical model category coincides with the realization model category when 
it exists.    

Denote the set of morphisms in
the homotopy category of the Reedy model category on $s\C$ by $[X,Y]^{\HR}$.
Call a map in $s\C$ a {\em realization weak equivalence} if
for all $Z$ in $\C$ it induces an isomorphism on $[-, cZ]^{\HR}$,
where $c$ is the constant functor.
An object in $s\C$ is {\em homotopically constant} if each of the simplicial
operators $d_i, s_i$ is a weak equivalence.  

\begin{theorem}\label{thm-A}
Let $\C$ be a model category.  Then there is at most one model
category structure on $s\C$ such that
\begin{itemize}
\item every level equivalence is a weak equivalence,
\item the cofibrations are the Reedy cofibrations, and
\item the fibrant objects are the homotopically constant, 
Reedy fibrant objects. 
\end{itemize}
When this model category exists, we refer to it as the {\em canonical
model category on $s\C$.}  Moreover, when it exists the weak equivalences
coincide with the realization weak equivalences.
\end{theorem}

\begin{proof}
First assume this canonical model category exists. Then since Reedy 
cofibrations are cofibrations and level equivalences are weak equivalences, a
Reedy cylinder object (\cite[I.1 Def.\ 4]{Q}, \cite[1.2.4]{hovey}) 
for a Reedy cofibrant object is also a cylinder
object in the canonical model category.  
This shows using~\cite[I.1 Cor.\ 1]{Q} that for $A$ Reedy cofibrant 
and $X$ homotopically constant and Reedy fibrant 
the homotopy classes of maps coincide in the homotopy category 
of the Reedy model category and the homotopy category of the 
canonical model category, $[A,X]^{\HR} \iso [A,X]^{\HC}$. 
Since level equivalences are weak equivalences
in both cases this means that for arbitrary $A$ and homotopically constant
$X$, $[A,X]^{\HR} \iso [A,X]^{\HC}$.  

A map $f \s A \to B$ is a weak equivalence in the canonical model category
if and only if for each homotopically constant $X$,
$[f, X]^{\HC}$ is a bijection.  Or, equivalently, $[f,X]^{\HR}$ is
a bijection.  Since $X$ is level equivalent to $c(X_0)$, this is equivalent
to $[f,cZ]^{\HR}$ being a bijection for each $Z$ in $\C$.  So the weak
equivalences are the realization weak equivalences.  

Since the cofibrations and weak equivalences are determined, the fibrations 
are determined by the right lifting property.  Hence there is at most
one model category on $s\C$ with the above properties.
\end{proof}

This specifies the model category of interest on $s\C$ because
when the canonical model category exists on $s\C$ it is Quillen 
equivalent to the original model category $\C$, see 
Proposition~\ref{prop-nice}.  

\begin{remark}\label{rem-hocolim}
In~\cite[21.1]{CS} and~\cite[21]{HH}, for any model category $\C$ a homotopy 
colimit functor is constructed which is the total left derived
functor of colimit.   Using this definition we could have defined the 
realization weak equivalences as those maps whose 
homotopy colimit is an isomorphism.    We use ``realization" instead of 
``hocolim" to avoid conflict with the terminology of~\cite{dugger}.  
Specifically, let 
$\hocolim \mc \HR \to \Ho(\C)$ be the total left derived functor of 
colimit.  Then $[A, cZ]^{\HR}$ is isomorphic to $[\hocolim A, Z]^{\Ho(\C)}.$ 
So $f\mc A \to B$ is a realization weak equivalence if and only if
$\hocolim f$ is an isomorphism.
In the rest of this paper though we only assume the existence of
the homotopy colimit for simplicial model categories, which follows 
from~\cite[XII]{BK}, see also~\cite[20]{HH}.  
\end{remark}

Now we demonstrate conditions which ensure the existence of the canonical model
category structure on $s\C$.  

\begin{definition}
A Reedy fibration $f \mc X \to Y$ in $s\C$ is an 
{\em equifibered Reedy fibration}
if the map $X_{m+1} \xrightarrow{(d_i,f_{m+1})} X_m \times_{Y_m} Y_{m+1}$ 
is a weak equivalence for each $m$ and for each simplicial face operator
$d_i$ with $0 \leq i \leq m+1$.
\end{definition}

\begin{axiom}[Realization Axiom]\label{real-axiom}
If $f \mc X \to Y$ in $s\C$ is an equifibered Reedy fibration and
a realization weak equivalence then $f$ is a level weak equivalence.
\end{axiom}

See Section~\ref{examples} for examples
where Axiom~\ref{real-axiom} is verified.
See Definition~\ref{def-cof-gen} for a definition of a cofibrantly generated
model category.  A model category is
{\em left proper} if the pushout of a weak equivalence along a cofibration
is a weak equivalence.  
Let $\ev \mc s\C \to \C$ be the evaluation functor 
given by $\ev X = X_0$.  Note $\ev$ is right adjoint to $c$, the constant
functor.

\begin{definition}  A pair $L, R$ of adjoint functors between
two model categories is a {\em Quillen adjoint pair} if $L$, the
left adjoint, preserves cofibrations and trivial cofibrations.  Equivalently,
$R$ preserves fibrations and trivial fibrations.   Such an adjoint pair
induces adjoint total derived functors on the homotopy categories, 
see~\cite[I.4 Thm.\ 3]{Q}.  A Quillen adjoint pair is a {\em Quillen
equivalence} if the total derived functors induce an equivalence on the
homotopy categories.  
\end{definition}

\begin{theorem}\label{main-thm}
If $\C$ is a left proper, cofibrantly generated model category that 
satisfies the Realization Axiom, then the following hold.
\begin{enumerate}
\item\label{thm-simplicial-model} The canonical model category on $s\C$ exists.
Moreover, it is cofibrantly generated and the fibrations are the
equifibered Reedy fibrations.  It is also   
referred to as the {\em \simp model category}.
\item \label{thm-simplicial} 
The \simp model category structure on $s\C$ 
satisfies Axiom~\ref{SM7} (SM7).  Hence it is a simplicial model category.
\item \label{thm-Quillen}
The adjoint functor pair $c \mc \C \leftrightarrows s\C \mc \ev$ induces a 
Quillen equivalence of the model category on $\C$ and the \simp model
category on $s\C$.  
\end{enumerate}
Moreover, the \simp model category structure agrees with the canonical
model category on $s\C$. 
\end{theorem}

This theorem is proved in Section~\ref{proofs}.  Recall that our definition
of cofibrantly generated is slightly stronger than standard; see
Definition~\ref{def-cof-gen}.  Since the weak
equivalences and cofibrations of the realization model category agree with
those of the canonical model category, these two model categories agree when 
they exist.  Thus, under the hypotheses of this theorem, the canonical model 
category is a simplicial model category.  
In fact, one can show that if the canonical
model category exists and is cofibrantly generated in the sense of 
Definition~\ref{def-cof-gen} then it is a simplicial model category. 

The next proposition shows that Realization Axiom~\ref{real-axiom} must hold 
if the fibrations in the canonical model category on $s\C$ are 
to be the equifibered Reedy fibrations.  

\begin{proposition}\label{prop-fib-real}
Assume $\C$ is a left proper, cofibrantly generated model category and the 
canonical model category on $s\C$ exists. 
Then the fibrations in the canonical model structure 
coincide with the equifibered Reedy fibrations 
if and only if $\C$ satisfies Realization Axiom~\ref{real-axiom}.
\end{proposition}

\begin{proof}
If the Realization Axiom holds, then part 1 of Theorem~\ref{main-thm}
gives the characterization of the fibrations as equifibered Reedy fibrations.
For the other implication, an equifibered Reedy fibration 
that is also a realization weak equivalence is a trivial fibration 
in the canonical model structure by assumption.  
But a trivial fibration has the right lifting property
with respect to the Reedy cofibrations, and hence is a level equivalence.  
Thus the Realization Axiom holds.
\end{proof}

\begin{remark}\label{rem-dugger}
As mentioned in the introduction, Dugger~\cite{dugger} also has conditions on a 
model category $\C$ which ensure that $s\C$ has a model
category structure, called the hocolim model category, which agrees with
the canonical model category and is simplicial.  
In particular, Proposition~\ref{prop-fib-real} can be used to 
explicitly describe the fibrations for some of Dugger's examples.  
\end{remark}

We end this section by stating a few of the properties that follow just from
the existence of the canonical model category structure.
Note that Theorem~\ref{main-thm}~(\ref{thm-Quillen}) follows from 
Theorem~\ref{main-thm}~(\ref{thm-simplicial-model}) 
and the first statement below since the realization
model category and the canonical model category agree when they exist.

\begin{proposition}\label{prop-nice}
If the canonical model category on $s\C$ exists then
\begin{enumerate}
\item The model category on $\C$ is Quillen equivalent to the canonical 
model category on $s\C$ via the adjoint functor pair $(c, \ev)$.
\item
A map between fibrant objects is a weak equivalence if and only if it is a 
level equivalence.
\item
The fibrations between fibrant objects are the Reedy fibrations. 
\end{enumerate}
\end{proposition}

\begin{proof}
For the second statement, note that $c$ preserves cofibrations and trivial cofibrations. 
By adjointness $\ev$ preserves fibrations and trivial fibrations, and hence also weak 
equivalences between fibrant objects.  
But, if $\ev f$ is a weak equivalence then $f$ is a level equivalence
since fibrant objects are homotopically constant.  

To show that the adjoint functor pair $(c, \ev)$ induces a Quillen equivalence,
we use the criterion in~\cite[4.1.7]{hss} since $\ev$
preserves and detects weak equivalences between fibrant objects.  
So we must show for any cofibrant object $X$ in $\C$ that $X \to \ev (cX)^f$ 
is a weak equivalence where $(cX)^f$ is a fibrant replacement of $cX$ in $s\C$.
Take $(cX)^f$ to be 
the Reedy fibrant replacement of $cX$, it is homotopically constant and
hence also a fibrant replacement in the canonical model category.  
Then $(cX)^f$  and $cX$ are level equivalent so $X \to \ev (cX)^f$ is indeed a 
weak equivalence in $\C$.

Since fibrations have the right lifting property with respect to level trivial Reedy 
cofibrations, a fibration is a Reedy fibration.  So 
we assume $f\mc X \to Y$ is a Reedy fibration between two fibrant objects and 
show that it is a fibration.   Factor $f= pi$ 
with $i$ a trivial cofibration and $p$ a fibration.
Then $i$ is a weak equivalence between fibrant objects, hence a level 
equivalence by part two.  Thus $i$ is a trivial Reedy cofibration
so it has the left lifting property with respect to $f$.  This implies that
$f$ is a retract of $p$, and hence a fibration in $s\C$. 
\end{proof}

\section{Examples}\label{examples}
In this section we give a criterion for simplicial model categories to
satisfy the Realization Axiom and verify the \real Axiom for stable model 
categories. So for the left proper, cofibrantly generated model categories 
among these examples, Theorem~\ref{main-thm} shows that $\C$ is 
Quillen equivalent to the simplicial, canonical model category on $s\C$.    
We mention one particular example, the category $\D$ of
unbounded differential graded modules over a differential graded
algebra.  

\subsection*{Simplicial model categories}
One source of model categories satisfying Realization Axiom~\ref{real-axiom}
is given by simplicial model categories where the realization factors through
simplicial sets, see below.   These examples are of interest for 
Sections~\ref{sec-unique} and~\ref{sec-functors},
where we discuss replacing functors between simplicial model categories
by simplicial functors and discuss
the uniqueness of simplicial model category structures.

For a simplicial model category $\C$, define a functor $\Sing \s \C \to s\C$
by $(\Sing X)_n = X^{\Delta[n]}$.  
Then $|-| \s s\C \to \C$ is the left adjoint to $\Sing$.
These functors are investigated further in Section~\ref{sec-unique}. 

\begin{definition}\label{def-factors}
For a simplicial model category $\C$, say that {\em the realization
factors through simplicial sets} if the following hold. 
\begin{enumerate}
\item There is a functor $U \s \C \to \Ss$ such that $f$ is a weak equivalence
in $\C$ if and only if $Uf$ is a weak equivalence in $\Ss$.
\item $U$ preserves fibrations. 
\item For any object $X \in s\C$,  $U|X|$ is naturally weakly equivalent to
$|\bU X|$ where $\bU$ is the prolongation of $U$ defined by applying $U$ to
each level in $s\C$. 
\end{enumerate}
\end{definition}

Examples of such model categories include topological spaces with 
$U= \Sing$ and the standard model category on simplicial objects in a category
$\C$ with an underlying set functor, such as simplicial groups~\cite[II.4]{Q}.  

A model category is {\em right proper} if
the pullback of a weak equivalence along a fibration is a weak equivalence.
A {\em proper} model category is one that is both right and left proper.

\begin{proposition}\label{prop-simp-real}
If $\C$ is a proper, cofibrantly generated simplicial model category where the 
realization factors through simplicial sets, as above, then $\C$ satisfies 
Realization Axiom~\ref{real-axiom}.  Hence the canonical model category on
$s\C$ exists, is simplicial, and is Quillen equivalent to $\C$ by 
Theorem~\ref{main-thm}.
\end{proposition}

Hence, under these hypotheses on $\C$, the applications
in Sections~\ref{sec-unique} and~\ref{sec-functors} apply.
These statements basically follow because the Realization Axiom holds for 
simplicial sets.

\begin{lemma}\label{lem-simp-real}
The model category of simplicial sets, $\Ss$, satisfies Realization
Axiom~\ref{real-axiom}.
\end{lemma}

Below we verify that Lemma~\ref{lem-simp-real} is a special case of the 
following proposition, 
essentially due to Puppe \cite{puppe-remark-homotopy-fibrations}.

\begin{proposition}\label{prop-equifibered-spaces}
Let $I$ be a small category and $X\rightarrow Y$ be a map of
$I$-diagrams of simplicial sets
such that for each 
$i_1\rightarrow i_2 \in I$ the square
\[\begin{CD}
{X(i_1)} @>>>   {Y(i_1)} \\
@VVV @VVV\\
{X(i_2)} @>>>  {Y(i_2)}
\end{CD}
\]
is homotopy cartesian. Then for each object $i\in I$, the square
\[\begin{CD}
{X(i)} @>>>   {Y(i)} \\
@VVV @VVV\\
{\hocolim_I X } @>>>  {\hocolim_I Y}
\end{CD}
\]
is homotopy cartesian. 
\end{proposition} 

\begin{proof}[Proof of Lemma~\ref{lem-simp-real}]
In the proposition take $I=\Delta$, the simplicial indexing category.  
An equifibered Reedy fibration $f\colon X\rightarrow Y$,
viewed as a map of $\Delta$-diagrams, satisfies the hypotheses of 
Proposition \ref{prop-equifibered-spaces}, 
and $f$ is a realization weak equivalence precisely when 
$\hocolim_{\Delta} X\rightarrow \hocolim_{\Delta} Y$ is a weak equivalence
by Remark~\ref{rem-hocolim}.
Therefore, for such $f$ and for every $i\in \Delta$ the map
$X(i)\to Y(i)$ is a weak equivalence, i.e., $f$ is a level weak equivalence.
\end{proof}

A proof of Proposition \ref{prop-equifibered-spaces} in this generality
appears in \cite{rezk}
where it is generalized to simplicial sheaves.
Alternatively, one can adapt the argument of \cite[App. HL]{farjoun},
where the Proposition is stated under the additional
hypothesis that the nerve of the indexing category $I$ and
all $Y(i)$ are path-connected. This implies that the homotopy colimit
of $Y$ is also connected, and so the conclusion as given in
\cite[App. HL]{farjoun} in terms of homotopy fibres is equivalent to
the conclusion of Proposition \ref{prop-equifibered-spaces}.
Proposition \ref{prop-equifibered-spaces} avoids explicit reference 
to homotopy fibres, and in this form the connectivity hypotheses 
are irrelevant.
It can be proved, as in \cite[App. HL]{farjoun}, 
by first checking the special cases of a homotopy pushout, 
a (possibly infinite) disjoint union and a sequential homotopy colimit;
an arbitrary homotopy colimit is built from these three ingredients, so the
result follows.

Puppe's original result is about simplicial objects in the category of
topological spaces; we could have derived the Realization Axiom for simplicial 
sets directly from his result, although some care would be needed, since 
he effectively works in a different model category (in which the ``weak
equivalences'' of spaces are plain homotopy equivalences) and he uses
the version of geometric realization of simplicial spaces in which
degeneracies are not collapsed.

\begin{proof}[Proof of Proposition~\ref{prop-simp-real}]
Let $f \s X \to Y$ be an equifibered Reedy fibration and a realization weak 
equivalence in $s\C$.  Since $\C$ is a right proper model category, the
condition for an equifibered Reedy fibration is invariant under level
equivalences.  By definition level 
equivalences are realization equivalences.  Hence, we can assume that 
$X$ and $Y$ are Reedy cofibrant.  For simplicial model categories,
the realization, $|-|$ is weakly equivalent to the homotopy colimit on Reedy 
cofibrant objects.  This follows from the generalization of~\cite[XII]{BK} to
general simplicial model categories, see~\cite[20.6.1]{HH}.
So $|f|$ is a weak equivalence in $\C$ by Remark~\ref{rem-hocolim}, 
since $f$ is a realization weak equivalence.  
By properties (1) and (2) of Definition~\ref{def-factors}, this means that 
$U|f|$ and $|\bU f|$
are weak equivalences.  Thus,  $\bU f$ is a realization 
weak equivalence of bisimplicial sets, by Remark~\ref{rem-hocolim} and
the fact that all bisimplicial sets are Reedy cofibrant.  Since $U$ preserves 
fibrations
and weak equivalences, it preserves homotopy pullback squares, and hence 
$\bU$ preserves equifibered Reedy fibrations.  
So, by Lemma~\ref{lem-simp-real}, $\bU f$ is a level 
equivalence.  Thus $f$ is a level equivalence.
\end{proof}

\subsection*{Stable model categories}
Recall from ~\cite[I.2]{Q} that the homotopy category of a pointed
model category supports a suspension functor $\Sigma$ with a
right adjoint loop functor $\Omega$.
A pointed model category $\C$ is {\em stable} if 
$\Sigma$ and $\Omega$ are inverse equivalences on the homotopy
category.  

\begin{proposition}\label{prop-stable-axiom}
Any proper, cofibrantly generated, stable model category $\C$ satisfies 
\real Axiom~\ref{real-axiom}.  Hence the canonical model category on
$s\C$ exists, is simplicial, and is Quillen equivalent to $\C$ by 
Theorem~\ref{main-thm}.
\end{proposition}

\begin{proof}
First note that since $\C$ is stable 
the Reedy model category on $s\C$ is also stable.
This follows since Reedy cofibrations and fibrations are level cofibrations and
fibrations and colimits and limits are taken levelwise. 
So the suspension and loop functors in the Reedy model category are level
equivalent to the levelwise suspension and loop in $\C$.

Now given a \simp weak equivalence $f\mc X \to Y$ in $s\C$ that is an 
equifibered Reedy fibration, we must show that $f$ is a level equivalence.   
Since $\C$ is right proper, the level homotopy fiber of $f$ is weakly 
equivalent to $F$, the fiber of $f$.  
In a stable model category fiber sequences induce long exact sequences after
applying $[-,cZ]^{\HR}$.  So $[F, cZ]^{\HR}$ is trivial for any $Z$ in $\C$. 
Since $f$ is equifibered, $F$ is homotopically constant and hence level 
equivalent to $c(F_0)$.  Thus $\id_F$ is trivial in $\HR$.  This implies that 
$F$ is level trivial, and hence that $f$ is a level equivalence since 
$\C$ is stable.
\end{proof}

\subsection*{Differential graded modules}
A cofibrantly generated model category, $\D$, of differential graded modules 
over a differential
graded algebra, $A$, is constructed in~\cite[5]{SS1}, see 
also~\cite[2.3.11]{hovey}.
The weak equivalences and fibrations are the quasi-isomorphisms and surjections
of the underlying $\Z$-graded chain complexes.
Since $\D$ is stable and proper, the realization axiom
follows by Proposition~\ref{prop-stable-axiom}.  
Thus, the following corollary follows from Theorem~\ref{main-thm}. 

\begin{corollary}\label{cor-chain complexes}
The proper, cofibrantly generated model category $\D$ of differential 
graded modules
over a differential graded algebra $A$ is Quillen equivalent to the
simplicial model category $s\D$ with the realization model category structure.
\end{corollary}

This answers a problem stated by Hovey,~\cite[8.9]{hovey}, 
which asks for a simple simplicial model category that is Quillen equivalent 
to unbounded chain complexes of $R$-modules, Ch($R$).  
Here $A$ is the differential graded algebra that is $R$ concentrated 
in degree zero. 

To make this example even more explicit, one can show that the {\em total
complex} functor $T$ is weakly equivalent to the homotopy colimit.
Let $X\in \sD$ be a simplicial object of differential graded $A$-modules.
We denote by $X_{s,t}$ the group in simplicial level $s$ and chain degree $t$. 
The total complex of $X$ is the chain complex
with levels $TX_n= \osum_{s+t=n} X_{s,t}$ and with total differential 
$d_{tot}=(-1)^s d + d'$.  Here $d$ is the internal chain differential 
in each simplicial level and $d'=\Sigma (-1)^i d_i$.
$TX$ is again a differential graded $A$-module.
Then a map $f$ is a realization weak equivalence in $s\D$ if and only if
$Tf$ is a quasi-isomorphism.

\section{Uniqueness of simplicial model category structures}\label{sec-unique}
In this section we consider categories $\C$ that already have
a given simplicial model category structure.  We then show that $\C$ is 
Quillen equivalent to $s\C$ via simplicial functors, see 
Theorem~\ref{thm-simp-Qui}.  
As stated in Corollary~\ref{cor-simp-Qui}, this 
implies that simplicial model category structures on a fixed model category are
unique up to simplicial Quillen equivalence.  See also
Corollary~\ref{cor-SMC-Qui} for a generalization of this result.
For these two statements we only need to assume that the canonical
model category on $s\C$ exists and is a simplicial model category.  We refer
to this as assuming the existence of the simplicial, canonical model category.
So the hypotheses considered in~\cite{dugger} work equally as well as the 
hypotheses considered in Theorem~\ref{main-thm}.  
Also, Proposition~\ref{prop-simp-real} provides many examples of simplicial 
model categories where the simplicial, canonical model category 
on $s\C$ exists.

First we recall the definition of a simplicial functor. 

\begin{definition}
Let $\C$ and $\D$ be categories enriched over simplicial sets.
Then a simplicial functor  $F \mc \C \to \D$ consists of a map
$F\mc \text{Ob}\,\C\to\text{Ob}\,\D$ of objects together with 
maps of simplicial sets $F\mc \map_{\C}(X,Y) \to \map_{\D}(FX, FY)$ 
that are associative and unital, see~\cite[II 1]{Q}.
\end{definition}

Since the vertices of the simplicial set $\map_{\C}(X,Y)$ are the morphisms
in the category $\C$, the restriction of a simplicial functor $F$ 
to vertices is an ordinary functor.
If the categories  $\C$ and $\D$ are also tensored over simplicial sets, 
then endowing an ordinary functor with a simplicial structure is equivalent 
to giving a transformation $K \otimes FX \to F(K \otimes X)$ 
that is natural in the simplicial set $K$ and in $X\in\C$ 
and that satisfies certain associativity and unity conditions, 
see~\cite[11.6]{HH}. 

For $\C$ a simplicial model category we now recall the adjoint functors
$\Sing \mc \C \to s\C$ and $|-| \mc s\C \to \C$.  For $X$ an object in $\C$, 
$\Sing(X)$ is the
simplicial object with $\Sing(X)_n= X^{\Delta[n]}$.  For $Y$ an object in $s\C$,
$|Y|$ is a coend~\cite[IX.6]{maclane} or 
the coequalizer of the following diagram induced by the simplicial
operators.  
\[
\amalg_{m,n}  X_m \otimes \Delta[n] \parallelarrows{1cm} \amalg_n  
X_n \otimes \Delta[n] 
\]
Throughout this section $X^K$, for $X$ in $\C$ and $K$ a simplicial set,
refers to the adjoint of the simplicial action on $\C$.  The simplicial
structure on $s\C$ is still as in Section~\ref{begin} and~\cite[II 1]{Q}.

\begin{theorem}\label{thm-simp-Qui}
Let $\C$ be a simplicial model category such that the simplicial, canonical
model category on $s\C$ exists.
Then the adjoint functors $\Sing$
and $|-|$ are simplicial and induce a Quillen equivalence between $\C$ and 
the simplicial, canonical model category structure on $s\C$.
\end{theorem}

Since the structures on $s\C$ are independent of any simplicial structure
on $\C$, this gives the following uniqueness statement for simplicial
model category structures.

\begin{corollary}\label{cor-simp-Qui}
Let $\C_1$ and $\C_2$ be two simplicial model categories with the
same underlying model category $\C$ such that the simplicial, canonical 
model category on $s\C$ exists.
Then $\C_1$ and $\C_2$ are simplicially Quillen equivalent.  
\end{corollary}

\begin{proof}
Apply Theorem~\ref{thm-simp-Qui} to both $\C_1$ and $\C_2$.  Then
they are both simplicially Quillen equivalent to $s\C$.
\end{proof}

To prove Theorem~\ref{thm-simp-Qui} we first prove that $\Sing$ and
$|-|$ are simplicial.

\begin{proposition}\label{prop-sing-simp}
For $\C$ a simplicial model category, $\Sing\mc \C \to s\C$ and $|-| \mc
s\C \to \C$ are simplicial functors. 
\end{proposition}

\begin{proof}
To show that $|-|$ is a simplicial functor we show that $K \otimes_{\C} |X|$
is isomorphic to $|K \otimes_{s\C} X|$.  
Here $\otimes_{\C}$ and $\otimes_{s\C}$
are the simplicial actions in the respective categories.  
These are not to be confused with the coends, see~\cite{maclane},  
$\otimes_{\Delta}$ and $\otimes_{\Delta \times \Delta}$ which follow.  
Since the left adjoint $|-|$ is a strong simplicial functor, that is,
the natural transformation is an isomorphism, it follows that the
right adjoint $\Sing$ is also a simplicial functor.

Let  $\mbD \s \Delta \to \Ss$ be the functor such that $\mbD(n)= \Delta[n]$,
the simplicial $n$-simplex.  Then $|X|$ is isomorphic to the coend 
$X \otimes_{\Delta} \mbD$ and for any simplicial set $K$, 
$K \iso (K \otimes_{\Delta} \mbD)$.  Because $\otimes_{\C}$
commutes with colimits, $K \otimes_{\C} |X| \iso (K \otimes_{\Delta} \mbD)
\otimes_{\C} (X \otimes_{\Delta} \mbD) \iso (K\cdot X) \otimes_{\Delta
\times \Delta} \mbD \times \mbD$.   
Here $(K \cdot X)(m, n) = K_m \cdot X_n$.  
The functor $\mbD \times \mbD$ is the left Kan extension of $\mbD$
across the diagonal functor $\delta \s \Delta \to \Delta \times \Delta$.
So $(K\cdot X) \otimes_{\Delta \times \Delta} \mbD \times \mbD
\iso \delta^*(K\cdot X) \otimes_{\Delta} \mbD$.  But $\delta^*(K \cdot X)$
is the functor describing $K \otimes_{s\C} X$, so this gives an isomorphism
of the last step with $|K \otimes_{s\C} X|$.  This produces the required
isomorphism. 
\end{proof}

\begin{proof}[Proof of Theorem~\ref{thm-simp-Qui}] 
First note that $M_n (\Sing X) = X^{\dot{\Delta}[n]}$
where $\dot{\Delta}[n]$ denotes the boundary of the simplicial $n$-simplex.  
So if $f \s X \to Y$ is a Reedy (trivial) fibration then 
$\Sing X \to \Sing Y$ is a Reedy (trivial) fibration because the induced map 
$X_n \to M_n X \times_{M_n Y} Y_n$
is equivalent to the map $X^{\Delta[n]} \to X^{\dot{\Delta}[n]} \times_{Y^{
\dot{\Delta}[n]}} Y^{\Delta[n]}$ which is a 
(trivial) fibration by the adjoint 
form of SM7, see SM7(a)~\cite[II 2]{Q}.  
The trivial fibrations in $s\C$ are the Reedy trivial fibrations.
Since the fibrations in $s\C$ 
between fibrant objects are Reedy fibrations by Proposition~\ref{prop-nice}, this 
shows that $\Sing$ preserves trivial fibrations and fibrations between fibrant objects.
Hence, by~\cite[A.2]{dugger}, $\Sing$ also preserves fibrations. 
By adjointness, $|-|$ preserves cofibrations and trivial cofibrations.

Since $|-|$ preserves trivial cofibrations it preserves weak equivalences 
between cofibrant objects.  It also detects weak equivalences between
cofibrant objects by Remark~\ref{rem-hocolim} since $|-|$ is weakly
equivalent to the homotopy colimit on Reedy cofibrant objects, by~\cite[XII]{BK} and~\cite[20.6.1]{HH}.  Hence by the dual of the criterion for
Quillen equivalences in~\cite[4.1.7]{hss}, we only need to check that for
fibrant objects $X$ in $\C$, $|(\Sing X)^c| \to X$ is a weak equivalence where
$(\Sing X)^c \to \Sing X$ is a trivial fibration from a cofibrant object in 
$s\C$.  By the simplicial
model category structure on $s\C$, $\Sing X$ is homotopically constant. 
Since 
$(\Sing X)^c$ is level equivalent to $\Sing X$, it is also homotopically
constant.   

Consider the following commuting 
square
\[
\begin{CD}
|c (\Sing X)^c_0| @>>> |(\Sing X)^c| @>>> |(\Sing X)^c| \\
@VVV @VVV @VVV\\
|cX| @>>> |\Sing X| @>>> X 
\end{CD}
\]
The left vertical map is a weak equivalence since $|cY| \iso Y$. The top
map is a weak equivalence since $(\Sing X)^c$ is homotopically constant. 
Finally, the bottom composite is the identity map.  Hence the right hand map 
is a weak equivalence as required.
\end{proof}

\section{Simplicial functors}\label{sec-functors}
In this section we again consider categories $\C$ 
that already have a given simplicial model category structure.
Since we have simplicial replacements 
for model categories, we now consider simplicial replacements of functors.
We show that a functor that preserves weak equivalences between
fibrant objects can be replaced by a simplicial functor that is weakly 
equivalent to the given functor on fibrant objects.  
We also show that a Quillen adjoint pair between simplicial model
categories can be replaced by a simplicial Quillen adjoint pair.
Combined with Theorem~\ref{thm-simp-Qui} this shows that if two simplicial
model categories have Quillen equivalent underlying model categories 
then they are in fact simplicially
Quillen equivalent, see Corollary~\ref{cor-SMC-Qui}.

For a functor $F\mc \C \to \D$, let $\bar{F}\mc s\C \to s\D$ be the 
prolongation of $F$ defined by applying $F$ at each level.  

\begin{proposition}\label{prop-Qui-adjoint}
Let $\C$ and $\D$ be model categories for which 
the simplicial, canonical model structures on $s\C$ and $s\D$ exist. 
Let $L\mc\C\to\D$ and $R\mc\D\to\C$ be a Quillen adjoint pair of functors. 
Then $\bar{L}$ and $\bar{R}$ are a simplicial Quillen adjoint pair between the 
simplicial model categories $s\C$ and $s\D$.  Moreover, if $L, R$ form a 
Quillen equivalence, so do $\bar{L}, \bar{R}$.
\end{proposition}

This answers Hovey's question in~\cite[8.9]{hovey} about replacing
Quillen adjunctions by Quillen equivalent simplicial Quillen adjunctions.
Indeed, if  $\C$ and $\D$ are {\em simplicial} model categories,
then Theorem~\ref{thm-simp-Qui} shows that $\C$ and $\D$ 
are simplicially Quillen equivalent to $s\C$ and $s\D$. 
So using Proposition \ref{prop-Qui-adjoint} one can replace a 
Quillen adjunction by a zig-zag of simplicial Quillen adjunctions
through $s\C$ and $s\D$ where the ``backwards" adjunction is a Quillen
equivalence.

\begin{proof}
First $\bar{L}$ is a simplicial functor. The necessary natural transformation,
$\bar{L}(X) \otimes K \to \bar{L}(X \otimes K)$ is given on each level
by the canonical maps $\coprod_{\sigma \in K_n} L(X_n) \to 
L(\coprod_{\sigma \in K_n} X_n)$.  

Since $R$ preserves fibrations, trivial fibrations, and 
limits, $\bar{R}$ preserves Reedy fibrations and Reedy trivial
fibrations.  So $\bar{R}$ preserves trivial fibrations and
fibrations between fibrant objects. By~\cite[A.2]{dugger} this implies
$\bar{R}$ also preserves fibrations.
Hence $\bar{L}, \bar{R}$ are a Quillen adjoint pair.   The last statement
follows from Theorem~\ref{thm-simp-Qui} and the two out of three property
for equivalences of categories, since Quillen
equivalences are Quillen adjoint functors that induce equivalences of homotopy 
categories~\cite[1.3.13]{hovey}.
\end{proof}

\begin{corollary}\label{cor-SMC-Qui}
Suppose that $\C$ and $\D$ are simplicial model categories for which 
the simplicial, canonical model structures on $s\C$ and $s\D$ exist.
 If there is a Quillen equivalence between the underlying model
categories $\C$ and $\D$, then $\C$ and $\D$ are simplicially Quillen 
equivalent.
\end{corollary}

\begin{proof}
By Theorem~\ref{thm-simp-Qui}, $\C$ and $\D$ are simplicially Quillen 
equivalent respectively to $s\C$ and $s\D$.  By 
Proposition~\ref{prop-Qui-adjoint}, the Quillen
equivalence between $\C$ and $\D$ can be lifted to a simplicial Quillen
equivalence between $s\C$ and $s\D$.
\end{proof}

Next we turn to constructing simplicial functor replacements.
Constructing simplicial cofibrant and fibrant replacement functors
is independent of the rest of this paper, see also~\cite[I.C.11]{farjoun} 
or~\cite{HH}. 
This construction is delayed to the end of the section.
These simplicial replacement functors are
then building blocks for replacing general functors by simplicial ones. 
In this section one can use the usual definition of cofibrantly generated
(see e.g.\ \cite[2.1.17]{hovey}), 
which is weaker than Definition~\ref{def-cof-gen}.

\begin{proposition}\label{prop-cof-rep}
For $\C$ any simplicial, cofibrantly generated model category there is a 
simplicial functorial factorization of any map $f\mc X \to Y$ as 
a cofibration followed by a trivial fibration and as a trivial cofibration
followed by a fibration.  In particular, this produces simplicial
cofibrant and fibrant replacement functors.  
\end{proposition}

\begin{proposition}\label{prop-simp-rep}
Assume $\C$, $\D$ are cofibrantly generated, simplicial model categories
such that the simplicial, canonical model categories on $s\C$ and $s\D$
exist and are cofibrantly generated.
Let $F\mc \C \to \D$ be a functor that preserves weak equivalences between 
fibrant objects.  Then $G=|Q\bar{F} \Sing (-)|$ is a simplicial functor, where 
$Q$ is a simplicial cofibrant replacement functor in the
simplicial, canonical model category on $s\D$.  Moreover,
there is a zig-zag of natural transformations between $F$ and $G$ that 
induce weak equivalences on fibrant objects in $\C$.
\end{proposition}

\begin{corollary}\label{all-wk-eq}
Assume $\C$, $\D$ are as above. 
If $F$ preserves all weak equivalences then $H=|Q{\bar {F}}\Sing R(-)|$ is a 
simplicial functor where $Q$ and $R$ are simplicial 
cofibrant and fibrant replacement functors in $s\D$ and $\C$ respectively.  
Moreover, for any $X$, $FX$ and $HX$ are naturally weakly equivalent.
\end{corollary}

\begin{proof}[Proof of Proposition~\ref{prop-simp-rep}] 
$G$ is a simplicial functor because each of its composites
is simplicial by Propositions~\ref{prop-sing-simp}, \ref{prop-Qui-adjoint}, 
and~\ref{prop-cof-rep}.

The first step in the zig-zag between $F$ and $G$ uses the
natural transformation $c \to \Sing$.  This induces 
$|Q\bar{F}c(-)| \to |Q\bar{F}\Sing (-)|=G(-).$ 
Note that for $X$ fibrant $cX \to \Sing X$ is a level equivalence 
between level fibrant objects by the simplicial model category structure on 
$\C$.  Since $|-|$ preserves trivial cofibrations by 
Theorem~\ref{thm-simp-Qui}, $|-|$ preserves weak equivalences between
cofibrant objects.
So, since $F$ preserves weak equivalences between fibrant objects, 
$|Q\bar{F}(cX)| \to |Q\bar{F}\Sing X|=GX$ is an equivalence for $X$ fibrant. 

To relate this to $FX$, note that $\bar{F}(cX)=cFX$.  
Since $QY \xrightarrow{p} Y$ is a level equivalence, 
$QcFX$ is homotopically constant.
Thus, $c\ev_0QcFX \to QcF$ is a level equivalence between cofibrant objects.
Hence $|c\ev_0QcFX| \to |Q\bar{F}cX|$ is also a weak equivalence for any $X$. 
$|c\ev_0QcFX| \to \ev_0QcFX$ is an isomorphism.  Since
$p$ is a level equivalence, $\ev_0QcFX \to FX$ is also an equivalence.
Combining this with the first step finishes the proof. 
\end{proof}

\begin{proof}[Proof of Proposition~\ref{prop-cof-rep}]
Given $f \mc X \to Y$ in $\C$
we construct a simplicial functorial factorization, $X \to Ff \to Y$, as a 
cofibration followed by a trivial fibration.  The other factorization is
similar.  Let $I$ be a set of generating cofibrations in $\C$.
Define the first stage, $F^1f$, as the pushout in the
following square.
\[
\begin{CD}
{ \amalg_{A_i \to B_i \in I} A_i \otimes 
(\map_{\C}(A_i, X)} \times_{\map_{\C}(A_i, Y)} \map_{\C}(B_i, Y)) @>>> X\\
@VVV @VVV\\
{ \amalg_{A_i \to B_i \in I} B_i \otimes 
(\map_{\C}(A_i, X)} \times_{\map_{\C}(A_i, Y)} \map_{\C}(B_i, Y)) @>>> F^1f\\
\end{CD}
\]
By~\cite[12.4.23]{HH}, any object that is small with respect to
the regular $I$-cofibrations is small with respect to all cofibrations.
So each $A_i$ is small relative to the cofibrations. 
Let $\kappa$ be the regular cardinal such that each $A_i$ is $\kappa$-small
with respect to the cofibrations.
Let $F^{\alpha + 1}f=F^1 (F^{\alpha}f \to Y)$ and for any limit ordinal $\beta <
\kappa$ let  $F^{\beta}=\colim_{\alpha} F^{\alpha}$.   Then we claim that
$F= F^{\kappa}$ is a cofibrant replacement functor which is also a simplicial
functor.

We need to show that $X \to Ff$
is a cofibration and that $Ff \to Y$ is a trivial fibration.  Since $\C$ is
a simplicial model category the left map in the square above is a 
cofibration.
Since pushouts and colimits preserve cofibrations this shows that
$X \to Ff$ is a cofibration.  To show that $Ff \to Y$ is a trivial fibration
we need to show that it has the right lifting property with
respect to any map $A_i \to B_i \in I$. 
Because $A_i$ is $\kappa$-small with respect to cofibrations,  the map 
$A_i \to Ff$ factors through
some stage, $F^{\alpha}f$.  Then, by construction, there is a lift 
$B_i\to F^{\alpha +1}f \to Ff$.

We now show that $F$ is simplicial.  The colimit of a diagram of
simplicial functors is again a simplicial functor.
Since the composition of simplicial functors is again simplicial,
we only need to show that $F^1$ is a simplicial functor. 
But $F^1$ itself is a colimit of functors which are simplicial, so we are done.
\end{proof}

\section{Reedy model category}\label{sec-Reedy}
In this section we show that the Reedy model category satisfies conditions
(1) and (2) but not (3) of Axiom~\ref{SM7}, (SM7).  These properties 
are also used in the proofs in Section~\ref{proofs}. 

The simplicial structure defined at the beginning of Section~\ref{main},
as with any simplicial structure, can be extended to morphisms.  Using
this structure on morphisms simplifies some of the notation and adjointness
properties that come up in verifying Axiom~\ref{SM7}, (SM7), for both the
Reedy and \simp model categories. 
See~\cite[5.3]{hss} for more about this structure on morphisms.

\begin{definition}
Given $f \mc X \to Y \in s\C$ and $i \mc K \to L \in \Ss$ define the {\em 
pushout product} 
$f \boxprod i \mc X \otimes L \coprod_{X \otimes K} Y \otimes K \to
Y \otimes L$ as 
the natural map from the pushout.  For $f$ in $\C$ define $f\boxprod i$ as 
$cf \boxprod i$ where $c\mc\C \to s\C$ is the constant functor.  Define
$f^{\boxprod i} \mc X^L \to Y^L \prod_{Y^K} X^K$ as the natural map to the 
pullback in $s\C$ or its zeroth level in $\C$ where the 
context will determine which category is meant.  
\end{definition}

Note that using this definition the map that appears in Axiom~\ref{SM7}, (SM7),
can be rewritten as the pushout product, $q=f\boxprod i$.  Also, note
that $- \boxprod i$ is adjoint to $(-)^{i}$. 
Next we rewrite the matching maps using this new notation.
Since $X^{\Delta[n]}=X_n$ and $X^{\dot\Delta[n]}=M_nX$, we have 

\begin{lemma}\label{lem-match} 
Let $f\mc X \to Y$ be a map in $s\C$.  The matching map $M_n f\mc X_n
\to Y_n \times_{M_nY} M_n X$ is the map $f^{\boxprod i_n}$ 
with $i_n\mc \dot\Delta[n] \to\Delta[n]$ the boundary inclusion.
\end{lemma}

\begin{proposition}\label{prop-proof}
If $g$ is a Reedy (trivial) fibration and $i$ is a cofibration in $\Ss$
then $g^{\boxprod i}$  in $s\C$ is a Reedy (trivial) fibration and hence
its zeroth level $g^{\boxprod i}$ in $\C$ is a (trivial) fibration.  
\end{proposition}

\begin{proof}
We need to consider the matching maps of $g^{\boxprod i}$, that is
$(g^{\boxprod i})^{\boxprod i_n}$ in $\C$ by Lemma~\ref{lem-match}.
Since $i \boxprod i_n$ is a cofibration in $\Ss$, it is enough to
show that $g^{\boxprod i}$ is a (trivial) fibration in $\C$.  In
fact it is enough to show this for each $i_n$ since they
generate the cofibrations in $\Ss$ by~\cite[3.2.2]{hovey}.  
But $g^{\boxprod i_n}$ is
a (trivial) fibration by Lemma~\ref{lem-match} since $g$ is a Reedy 
(trivial) fibration.
\end{proof}

A corollary of this Proposition is 
that although the Reedy model category is not simplicial
it does satisfy the first two properties of Axiom~\ref{SM7}, (SM7).

\begin{corollary}\label{cor-Reedy-simp}
Given $f\mc X\to Y$ a Reedy cofibration in $s\C$ and $i \mc K \to L$ a
cofibration in $\Ss$ then $f \boxprod i \mc X \otimes L \coprod_{X\otimes K}
Y\otimes K \to Y\otimes L$ is a Reedy cofibration.  Moreover, if $f$ is
also a level weak equivalence, then so is $f \boxprod i$.  But if $i$ is a
weak equivalence and $f$ is not, then $f\boxprod i$ is not necessarily
a weak equivalence.
\end{corollary}

\begin{proof}
The first two statements follow by adjointness from 
Proposition~\ref{prop-proof}.   For all three statements, see 
also~\cite[2.6]{e2} and compare with~\cite[5.4.1]{hovey}. 
\end{proof}

\section{Realization model category}\label{proofs}
In this section we prove Theorem~\ref{main-thm}, which states that the
\simp model structure on $s\C$ is a model category that is simplicial
and Quillen equivalent to the original model category $\C$.  

To verify the axioms for the realization model category on $s\C$ we 
assume that $\C$ is a cofibrantly generated model
category.   We now recall a version of the definition of 
cofibrantly generated model category from~\cite{DHK}, or 
see~\cite[2.1.17]{hovey}, \cite[2.2]{SS1}, or~\cite{HH}.  For a cocomplete
category $\C$ and a class $I$ of maps, the {\em $I$-injectives} are
the maps with the right lifting property with respect to the maps in $I$.
The {\em $I$-cofibrations} are the  maps with the left lifting property
with respect to the $I$-injectives.  
Finally, the {\em regular $I$-cofibrations}
(called {\em relative $I$-cell complexes} in \cite[2.1]{hovey})
are the (possibly transfinite) compositions of pushouts of maps in $I$.
In particular all isomorphisms are regular $I$-cofibrations, 
see the remark following~\cite[2.1.9]{hovey}.

\begin{definition}\label{def-cof-gen}
A model category $\C$ is {\em cofibrantly generated} if it is complete
and cocomplete and there exists a set of cofibrations $I$ and a set
of trivial cofibrations $J$ such that
\begin{enumerate}
\item the fibrations are precisely the $J$-injectives,
\item the acyclic fibrations are precisely the $I$-injectives,
\item the domain and range of each map in $I$ 
and each map in $J$ is {\em small} relative to the regular 
$I$-cofibrations, and
\item the domain and range of each map in $I$ is cofibrant.
\end{enumerate}
Moreover, here the (trivial) cofibrations are the $I$ ($J$)-cofibrations. 
\end{definition}   

For the definition of {\em small} see the above mentioned references.
The crucial reason for requiring a cofibrantly generated model category is 
the small object argument, Proposition~\ref{small-objects}, as in~\cite{Q}, 
see also~\cite{DHK} or ~\cite[2.1.14]{hovey}.
The smallness requirements here are stronger than what is necessary
for the small object argument to apply to $I$ and $J$; 
we added the requirement that the ranges of $I$ and $J$ are also small.  
We use this to show that the domains
of the new generators defined in~\ref{def-J-simp} for $s\C$ have small domains
so the small object argument will apply in $s\C$.
Since $\C$ is also assumed to be left proper, we could replace $J$ by
a set $J'$ of regular $I$-cofibrations and the smallness
condition for $J'$ would follow by~\cite[12.3.8]{HH}.
The maps in $I$ 
are required to be between cofibrant objects so that Proposition~\ref{prop-SM7}
holds.

\begin{proposition}[Small object argument]
\label{small-objects}
Let $\C$ be a cocomplete category and $I$ a set of maps in $\C$ whose
domains are small relative to the regular $I$-cofibrations.  Then
\begin{enumerate}
\item 
there is a functorial factorization of any map $f$ in $\C$ as $f=pi$
with $p$ an $I$-injective and $i$ a regular $I$-cofibration. 
And thus,
\item every $I$-cofibration is a retract of a regular $I$-cofibration.
\end{enumerate}
\end{proposition}

We now begin to verify the model category axioms
for the \simp model structure on $s\C$.  We assume that $\C$ is a left proper, 
cofibrantly generated model category that satisfies the Realization 
Axiom~\ref{real-axiom}.  For the factorizations we use 
Proposition~\ref{small-objects}.  We characterize the (trivial) fibrations
as the maps with the right lifting property with respect to a set of maps, $J$ 
($I$).  Let $I_C$ be a
set of generating cofibrations for $\C$ and $J_C$ be a set of generating
trivial cofibrations for $\C$.
In the category of simplicial sets, let $I_{\del}$ be the 
set of inclusions of boundaries
into simplices, $i_n\mc \dot\Delta[n] \to \Delta[n]$ for each $n$.  Let
$I_F$ be the set of inclusions of faces into simplices, $\delta_i\mc
\Delta[m] \to \Delta[m+1]$ for each $m$ and $0\leq i \leq m+1$. 

\begin{definition}\label{def-J-simp}
Let $I=I_C \boxprod I_{\del}$ denote the set of maps 
\[
f\boxprod i_n\mc A \otimes \Delta[n] \coprod_{A \otimes \dot{\Delta}[n]} B \otimes \dot{\Delta}[n] \to B \otimes \Delta[n]
\]
for each $n$ and $f\mc A \to B$ any map in $I_C$.  
Let $J'=J_C \boxprod I_{\del}$ denote the set of maps 
\[
f\boxprod i_n\mc A \otimes \Delta[n] \coprod_{A \otimes \dot{\Delta}[n]} B \otimes \dot{\Delta}[n] \to B \otimes \Delta[n]
\]
for each $n$ and $f\mc A \to B$ any map in $J_C$. 
Let $J''=I_C \boxprod I_{F}$ denote the set of maps 
\[
f\boxprod {\delta}_i\mc A \otimes \Delta[m+1] \coprod_{A \otimes {\Delta
[m]}} B \otimes {\Delta[m]} \to B \otimes \Delta[m+1] 
\]
for each $m$ and $i$ 
with $0\leq i \leq m+1$ and 
$f\mc A \to B$ any map in 
$I_C$.  Let $J$ be the union of the two sets $J'$ and $J''$.
\end{definition}

\begin{lemma}\label{lem-J-small}
The domains of $I$ and $J$ are small relative to the regular $I$-cofibrations.
\end{lemma}

\begin{proof}
We prove the statement for $J$, the statement for $I$ follows similarly. 
A finite colimit of small objects is small since finite limits commute
with small filtered colimits,~\cite[IX 2]{maclane}.  
The domains of $J$ can be built by
finite colimits from objects $X \otimes \Delta[n]$ for $X$ a domain
or range of a map in $I_{\C}$ or $J_{\C}$.  
Since $s\C(X \otimes \Delta[n], Y) \iso \C(X, Y^{\Delta[n]})\iso\C(X, Y_n)$
and $X$ is small relative to regular $I_{\C}$-cofibrations by 
Definition~\ref{def-cof-gen}, 
$X \otimes \Delta[n]$ is small relative to maps in $s\C$ that are regular 
$I_{\C}$-cofibrations on each level.  But each level of
a regular $I$-cofibration is a regular 
$I_{\C}$-cofibration.  This is because each level of a map in $I$ 
is just a direct sum of copies of maps in $I_{\C}$ or identity maps.  
Identity maps and coproducts of regular cofibrations are regular cofibrations.
So each level of each map in $I$ is a regular $I_{\C}$-cofibration. 
Hence this is also true of the regular $I$-cofibrations. 
\end{proof}

Since the domains are small we can use the 
small object argument, Proposition~\ref{small-objects}, 
to factor any map into an $I$ ($J$)-cofibration followed by an 
$I$ ($J$)-injective.  This applies directly to $I$ by Lemma~\ref{lem-J-small}. 
For $J$, since the domains of $J$ are small relative to
the regular $I$-cofibrations, they are small with respect to all 
cofibrations including the regular $J$-cofibrations by~\cite[13.3.3]{HH}.  
Hence Proposition~\ref{small-objects} applies.
To see that this gives us the needed factorization we show in the next 
propositions that an $I$ ($J$)-cofibration is a realization (trivial) 
cofibration and that a $J$ ($I$)-injective is a realization (trivial) 
fibration.  

\begin{proposition}\label{prop-simp-fib}
The $J$-injective maps are the equifibered Reedy fibrations.  In other words, the 
equifibered Reedy fibrations are the maps with the right lifting
property with respect to $J$.  The Reedy fibrations
are the maps with the right lifting property with respect to $J'$.  
Moreover, the $J$-injective objects are the homotopically constant,
Reedy fibrant objects.
\end{proposition}

\begin{proof}
A  Reedy fibration is
a map $f$ whose matching maps are fibrations.  These matching maps are 
$f^{ \boxprod i_n}$ with $i_n \in I_{\del}$ by Lemma~\ref{lem-match}.  
That is, $f^{\boxprod i_n}$ has the right
lifting property with respect to each map in $J_C$.   By adjointness,
this is equivalent to $f$ having the right lifting property with respect
to the maps in $J_C \boxprod I_{\del} = J'$.
 
Given a Reedy fibration $f\mathcolon X \to Y$, then 
$f^{\boxprod \delta_i} \s X_{m +1}\to X_{m }\times_{Y_{m}}Y_{m+1}$ is a 
fibration by Proposition~\ref{prop-proof}.  
So a Reedy fibration $f$ is equifibered if and only
if $f^{\boxprod \delta_i}$ is a  trivial fibration. 
By adjunction $f^{\boxprod \delta_i}$ is a trivial fibration if and only if
$f$ has the right lifting property with respect to $J'' = I_{\C} \boxprod I_F$.
So $f$ is an equifibered Reedy fibration if and only if $f$  
has the right lifting property with respect to $J$.
The last statement of the proposition follows since $f \s Z \to *$ 
is an equifibered Reedy fibration if and only if $Z$ is Reedy fibrant
and for each $n$ and $i$ the map 
$d_i \s Z_{n+1} \to Z_n$ is a trivial fibration.  
\end{proof}

Next we turn to the $I$-cofibrations and $I$-injectives. 

\begin{proposition}\label{prop-I-cof}
The $I$-injective maps are the Reedy trivial fibrations.  
Also, the Reedy trivial fibrations are the equifibered Reedy fibrations 
that are also realization weak equivalences.  
Hence, the $I$-cofibrations are the Reedy cofibrations.
\end{proposition}

\begin{proof}
Much as in the previous proof, a map $f$ is a  Reedy trivial fibration 
if the matching maps $f^{ \boxprod i_n}$ are trivial fibrations.
That is $f^{\boxprod i_n}$ has the right
lifting property with respect to each map in $I_C$.   By adjointness,
this is equivalent to $f$ having the right lifting property with respect
to the maps in $I_C \boxprod I_{\del} = I$.   

By the Realization Axiom~\ref{real-axiom},
an equifibered Reedy fibration that is also a realization weak equivalence is
a level equivalence, and hence a Reedy trivial fibration.  Conversely, for
$f$ a Reedy trivial fibration, the maps $f_n \mc X_n \to
Y_n$ are trivial fibrations.  Since $f_{n+1}$ factors as $X_{n+1} 
\to X_n \times_{Y_n} Y_{n+1} \to Y_{n+1}$ and the second map here 
is the pull back of a trivial fibration, 
the map  $X_{n+1} \to X_n \times_{Y_n} Y_{n+1}$ is a weak equivalence.  
So a Reedy trivial fibration is equifibered.  
Then, since level equivalences are realization weak equivalences, 
this shows that a Reedy trivial fibration is a realization trivial fibration, 
i.e. an equifibered Reedy fibration 
that is also a realization weak equivalence.
\end{proof}

Now we are left with verifying that the $J$-cofibrations are Reedy cofibrations
and realization weak equivalences.

\begin{proposition}\label{prop-simp-cof}
A $J$-cofibration is a Reedy cofibration and a \simp weak equivalence.
\end{proposition}

\begin{proof}
A $J$-cofibration has the left lifting property with respect to the
$J$-injective maps, the equifibered Reedy fibrations.  Since any
Reedy fibration that is also a level equivalence is equifibered, a
$J$-cofibration has the left lifting property with respect to the
Reedy trivial fibrations.  Hence a $J$-cofibration is a Reedy cofibration.

Each $J$-cofibration is a retract of a directed colimit of pushouts of
maps in $J$ by Proposition~\ref{small-objects}.
The maps in $J'$ are level equivalences, hence the maps built from $J'$
are Reedy trivial cofibrations. These level equivalences
are \simp weak equivalences.  So we only need to consider 
$J''$-cofibrations.  Since the maps in $I_F$ are
trivial cofibrations of simplicial sets, they are $I_{\Lambda}$-cofibrations
where $I_{\Lambda}=\{\lambda_n \mc \Lambda^k[n] \to \Delta[n]\}$ is the
set of inclusions of the horns into simplices.  Hence $J''$-cofibrations
are $(I_{\C} \boxprod I_{\Lambda})$-cofibrations.
Below, in Proposition~\ref{prop-SM7}, we show that
any $(I_{\C} \boxprod I_{\Lambda})$-cofibration is a \simp weak equivalence.
\end{proof}

To finish our verification of the realization model category structure
we need to use a different characterization of the realization weak
equivalences.

\begin{definition}
A map $f'\mc A' \to B'$ is a {\em cofibrant replacement} of a map 
$f \mc A \to B$
if $A'$ and $B'$ are cofibrant objects, $f'$ is a cofibration, and there exist
level equivalences $i_A \mc A' \to A$ and $i_B \mc B' \to B$ such that
$f i_A = i_B f'$.
\end{definition}

\begin{proposition}\label{prop-map}
A map $f \mc A \to B$ in $s\C$ is a realization weak equivalence
if and only if for some cofibrant
replacement $f' \mc A' \to B'$, and
for each homotopically constant, Reedy fibrant object $Z$ in $s\C$,
$\map (B',Z) \to \map(A',Z)$ is a weak equivalence.
\end{proposition}

The following lemmas are used to prove this proposition.   

\begin{lemma}\label{lem-lambda}
The map $Z^{\lambda_n} \mc Z^{\Delta[n]} \to 
Z^{\Lambda^k[n]}$ in $\C$ is a trivial fibration for $Z$ any homotopically
constant Reedy fibrant object in $s\C$.      
\end{lemma}

\begin{proof}
$Z^{\lambda_n}$ is a fibration, by Corollary~\ref{cor-Reedy-simp}.
Since $\Lambda^k[1]=\Delta[0]$, $Z^{\lambda_1}$ is 
the map $d_k \mc Z_1 \to Z_0$, which is a trivial fibration since $Z$ is
a homotopically constant, Reedy fibrant object.  This proves 
the lemma for $n=1$.  We proceed by induction.

$ Z^{\Lambda^k[n]}$ is the pullback of a punctured $n$-cube where
each arrow is of the form $Z^{\delta_i} \s Z^{\Delta[m]} \to
Z^{\Delta[m-1]}$, that is, $Z_{m} \to Z_{m-1}$ for $m < n$.   These maps are 
fibrations by Corollary~\ref{cor-Reedy-simp} and they are weak equivalences 
because $Z$ is homotopically constant.  By induction the map from the object
at the puncture of each contained punctured $k$-cube, for $k < n$, to the  
pullback is a trivial fibration.  For any such punctured $n$-cube, the added
maps from the pullback are trivial fibrations.  That is, the maps
from the pullback, $Z^{\Lambda^k[n]}$, to each $ Z^{\Delta[n-1]}=Z_{n-1}$ 
are trivial fibrations.
Since each $\delta_i$ factors as $\Delta[n-1] \to \Lambda^k[n] \to \Delta[n]$,
this proves the lemma holds for $n$ by the two out of three property for weak 
equivalences.
\end{proof}

\begin{lemma}\label{lem-constant}
For $K$ any simplicial set and $Z$ any homotopically constant, Reedy fibrant 
object, $Z^K$ is homotopically constant and Reedy fibrant. 
\end{lemma}

\begin{proof}
First note that by an adjoint of SM7 (i), which is verified for
the Reedy model category in Corollary~\ref{cor-Reedy-simp}, $Z^K$ is 
Reedy fibrant.
Hence by Proposition~\ref{prop-simp-fib}, $Z^K$ is $J'$-injective 
and we only need to show that $Z^K$ is $J''$-injective to finish the proof. 

Here we say ``$(f,g)$ has the lifting property," as short hand for
$f$ has the left lifting property with respect to $g$.  This also
extends to sets of maps.  By Lemma~\ref{lem-lambda}, $(i, Z^{\lambda_n})$ has
the lifting property for $i$ in $I_{\C}$, $\lambda_n$ in $I_{\Lambda}$,
and $Z$ any homotopically constant, Reedy fibrant object.  Let $H$
be the class of maps $Z \to *$ for such $Z$.  Then,
by adjointness $(I_{\C} \boxprod I_{\Lambda}, H)$ has the lifting 
property.   But then pushouts, colimits and retracts of maps in 
$I_{\C}\boxprod I_{\Lambda}$ also have the left lifting property with
respect to $H$.  That is, $((I_{\C} \boxprod I_{\Lambda})$-cofibrations,
$H)$ has the lifting property.
For $i$ a cofibration and $j$ a trivial cofibration of simplicial sets , the 
pushout product $j \boxprod i$ is an $I_{\Lambda}$-cofibration.  
So $f \boxprod j \boxprod i$ is an 
$(I_{\C} \boxprod I_{\Lambda})$-cofibration for $f$ in $I_{\C}$.
Hence $(I_{\C} \boxprod I_{\Lambda} \boxprod I_{\del}, Z)$ has the lifting 
property.  
Consider the cofibration $i\s \emptyset \to K$.  By adjointness this shows that
$(I_{\C} \boxprod I_{\Lambda}, Z^K)$ has the lifting property. 
Hence $Z^K$ is $J''$-injective.
\end{proof}

\begin{proof}[Proof of Proposition~\ref{prop-map}]
Our first claim is that  $\pi_0 \map(A,X)$ is naturally isomorphic
to $[A, X]^{\HR}$ for $A$ Reedy cofibrant and $X$ homotopically constant 
and Reedy fibrant.
Indeed the maps $X\iso X^{\Delta[0]} \xrightarrow{f} X^{\Delta[1]}
\xrightarrow{p} X^{\Delta[0] \amalg \Delta[0]} \iso X \times X$
produce $X^{\Delta[1]}$ as a path object for $X$.  Here $f$ is a level
equivalence by Lemma~\ref{lem-constant} since it is a map between
homotopically constant objects whose zeroth level is given by the
equivalence $s_1\s X_0 \to X_1$ and Proposition~\ref{prop-proof} shows that
$p$ is a Reedy fibration. This implies the claim.

Since $f$ is a realization weak equivalence if and only if its cofibrant
replacement is, we can restrict
to the case when $f$ is its own cofibrant replacement.  Then
requiring that $\map(f,Z)$ is a weak
equivalence for all homotopically constant, Reedy fibrant objects $Z$ 
is equivalent to requiring that for all simplicial sets $K$,
$\pi_0 \map(K, \map(f,Z))\iso \pi_0\map(f,Z^K)$ is a bijection for all 
such $Z$. 
By Lemma~\ref{lem-constant} and the above, 
this is equivalent to $[B, Z^K]^{\HR} \to [A, Z^K]^{\HR}$
being a bijection for all such $K$ and $Z$.  

As $Z$ runs through all homotopically constant, Reedy fibrant objects and 
$K$ runs through 
all simplicial sets, $(Z^K)_0$ runs through all fibrant objects in $\C$.  
Since $c(Z^K)_0 \to Z^K$ is a level equivalence, this is equivalent to
$[B, cX]^{\HR} \to [A,cX]^{\HR}$ being a bijection for all $X$ in $\C$.
\end{proof}

The following proposition finishes the identification of 
the $J$-cofibrations as \simp weak equivalences.  It is also
useful in checking that $s\C$ is a simplicial model category.  
 
\begin{proposition}\label{prop-SM7}
Any $(I_{\C} \boxprod I_{\Lambda})$-cofibration is a \simp weak equivalence.
\end{proposition}

\begin{proof}
By the proof above of Lemma~\ref{lem-constant}, 
$(I_{\C} \boxprod I_{\Lambda} \boxprod I_{\del}, Z)$ has the lifting property 
for $Z$ homotopically constant and Reedy fibrant.    Then 
by adjointness, $(I_{\del}, \map(I_{\C} \boxprod I_{\Lambda},  Z))$ also has 
the lifting 
property for any such $Z$. That is, any map in 
$\map(I_{\C} \boxprod I_{\Lambda}, Z)$ is a trivial fibration.  
Since the maps in $I_{\C}$ are assumed to be between cofibrant objects,
the maps in $I_{\C} \boxprod I_{\Lambda}$ are Reedy cofibrations between Reedy 
cofibrant objects.  So they are their own cofibrant replacements.  
Hence the maps in $I_{\C} \boxprod I_{\Lambda}$
are realization weak equivalences by Proposition~\ref{prop-map}.   
Since the maps in $I_{\C} \boxprod I_{\Lambda}$
are Reedy cofibrations, to finish this proof it is enough to show that
Reedy cofibrations that are realization weak equivalences are preserved
under pushouts, directed colimits, and retracts.

Since $\C$ is left proper, if $g$ is a pushout of a Reedy cofibration $f$ then 
one can choose a cofibrant
replacement $g'$ for $g$ as a pushout of the cofibrant replacement $f'$
of $f$.   Hence $\map(g',Z)$ is a pullback of $\map(f',Z)$.  We show in
the next paragraph that if $f'$ is a Reedy cofibration then $\map(f',Z)$ is
a fibration.  So if $f$ is a Reedy cofibration and realization weak equivalence
then $\map(f',Z)$ and hence also $\map(g',Z)$ is a trivial fibration.  Thus, 
$g$ is a realization weak equivalence.
Since retracts and directed limits of trivial fibrations are also trivial 
fibrations,
it follows that retracts and directed colimits also preserve Reedy 
cofibrations that are realization weak equivalences. 

Since $(I_{\C} \boxprod I_{\del} \boxprod I_{\Lambda}, Z)$ has the 
lifting
property, so does $((I_{\C} \boxprod I_{\del})$-cofibrations, 
$Z^{I_\Lambda})$
for $Z$ any homotopically constant, Reedy fibrant object.
By adjointness this shows that for any Reedy cofibration $i$, $\map(i,Z)$
is a fibration since it has the right lifting property with respect to
$I_{\Lambda}$.  
\end{proof}

\begin{proof}[Proof of Theorem~\ref{main-thm}~(\ref{thm-simplicial-model})]
As always, we assume that $\C$ is a left proper, 
cofibrantly generated model category that satisfies Realization 
Axiom~\ref{real-axiom}.  The category $s\C$ has all
limits and colimits since $\C$ does.  The two out of three axiom for
weak equivalences and the retract axiom for the cofibrations and
weak equivalences are easily checked.  The retract axiom for fibrations
follows from Proposition~\ref{prop-simp-fib}.
The two factorizations follow from Propositions~\ref{prop-simp-fib},
~\ref{prop-I-cof} and~\ref{prop-simp-cof}
by Proposition~\ref{small-objects}.  One lifting property follows from 
Proposition~\ref{prop-I-cof} since the \simp trivial fibrations are the Reedy 
trivial fibrations.  So only the lifting of a \simp trivial
cofibration with respect to an equifibered Reedy fibration is left.  Assume
$f\mathcolon X \to Y$ is a Reedy cofibration and a \simp weak equivalence.
Factor $f=pi$ where $i$ is a $J$-cofibration and $p$ is $J$-injective.
Since $f$ and $i$ are \simp weak equivalences, $p$ is also a \simp weak
equivalence.  Since $f$ is a Reedy cofibration, Propositions~\ref{prop-simp-fib}
and~\ref{prop-I-cof} show that it has the left lifting property with respect 
to $p$.  Thus,  $f$
is a retract of $i$.  Hence $f$ is a $J$-cofibration and so it has the
left lifting property with respect to any equifibered Reedy fibration.  This 
finishes the proof that the \simp model structure on $s\C$ is a model category.
\end{proof}

\begin{corollary}\label{prop-simp-cof-gen}
Let $I$ and $J$ be as defined in Definition~\ref{def-J-simp}.
The \simp model category on $s\C$ is cofibrantly generated, with
$I$ a set of generating cofibrations and $J$ a set of generating trivial
cofibrations.
\end{corollary}

We now prove Theorem~\ref{main-thm}~(\ref{thm-simplicial}), which states  
that the \simp model category structure on $s\C$ satisfies Axiom~\ref{SM7}, 
(SM7).  Hence, it is a simplicial model category.  

\begin{proof}[Proof of Theorem~\ref{main-thm}~(\ref{thm-simplicial})]
Given $f\mc A \to B$ a Reedy cofibration in $s\C$ and
$i\mc K \to L$ a cofibration in $\Ss$, 
$f\boxprod i$ is a Reedy cofibration by Corollary~\ref{cor-Reedy-simp}.
So we are left with showing that if $f$ or $i$ is also a weak equivalence
then so is $f\boxprod i$.

First consider the case where $i$ is a trivial cofibration. 
Since the pushout product of a trivial cofibration and a cofibration
of simplicial sets is a trivial cofibration, 
$((I_{\C} \boxprod I_{\del})$-cofibrations) $\boxprod 
(I_{\Lambda}$-cofibrations) is contained
in $(I_{\C} \boxprod I_{\Lambda})$-cofibrations.
So by Proposition~\ref{prop-SM7}, $f \boxprod i$ is a realization weak 
equivalence for $f$ any Reedy cofibration.
 
Next consider the case where $f$ is a \simp weak equivalence.  Since trivial
cofibrations are preserved under pushouts, retracts and colimits, it is enough 
to show that for $f$ in $J$, $f \boxprod i$ is a realization weak equivalence.
For $f$ in $J'$ this follows from Corollary~\ref{cor-Reedy-simp}.  For
$f$ in $J''= I_{\C} \boxprod I_F$ this follows from the previous paragraph
by associativity, since the maps in $I_F$ are trivial cofibrations.
\end{proof}

Recall that the Quillen equivalence
of $\C$ and $s\C$, Theorem~\ref{main-thm} part~(\ref{thm-Quillen}), follows
from Proposition~\ref{prop-nice} since the realization model category agrees
with the canonical model category on $s\C$.

\ifx\undefined\bysame
\newcommand{\bysame}{\leavevmode\hbox to3em{\hrulefill}\,}
\fi

\end{document}